\documentclass[a4paper]{amsart}
\usepackage{amssymb}
\usepackage{algorithmic}
\usepackage{algorithm}
\usepackage[all]{xy}
\newtheorem{theorem}{Theorem}[section]

\newtheorem{lemma}[theorem]{Lemma}

\newtheorem{corollary}[theorem]{Corollary}
\theoremstyle{definition}

\newtheorem{definition}{Definition}[section]
\theoremstyle{remark}

\newtheorem*{ack}{Acknowledgment}

\newtheorem{example}{Example}

\newcommand{\field}[1]{\mathbb{#1}}
\newcommand{\bbbn}{\field{N}}
\newcommand{\bbbr}{\field{R}}

\newcommand{\vrightarrow}[1]{\
  \xy\xymatrix{*!{}\ar@{.>}[r]^{\scriptscriptstyle #1}&*!{}}\endxy\ }
\newcommand{\vlrightarrow}[2]{\
  \xy\xymatrix@+#1{*!{}\ar@{.>}[r]^{\scriptscriptstyle #2}&*!{}}\endxy\ }

\DeclareMathOperator\depth{{\rm td}}
\DeclareMathOperator\clos{{\rm clos}}
\DeclareMathOperator\height{{\rm height}}
\DeclareMathOperator\tw{{\rm tw}}

\newcommand{\rdens}[1]{\nabla_{#1}}

\newcommand{\card}[1]{\lvert{#1}\rvert}

\DeclareMathOperator\md{\Delta^{\rm --}}

\def\vG{\vec{G}}
\def\vH{\vec{H}}

\newcommand{\rdpath}[1]{\xy\xymatrix{*!{}\ar@{~2>}[r]^{#1}&*!{}}\endxy}
\newcommand{\ldpath}[1]{\xy\xymatrix{*!{}&*!{}\ar@{~2>}[l]_{#1}}\endxy}

\begin{document}

\title{Grad  and Classes with Bounded Expansion II. Algorithmic Aspects.}
\author{Jaroslav Ne\v
  set\v ril}
\thanks{Supported by grant 1M0021620808 of the Czech Ministry of Education}
\address{Department of Applied Mathematics\\
  and\\
  Institute of Theoretical Computer Science (ITI)\\
  Charles University\\
  Malostransk\' e n\' am.25, 11800 Praha 1\\
  Czech Republic} \email{nesetril@kam.ms.mff.cuni.cz} \author{Patrice
  Ossona de Mendez}
\address{Centre d'Analyse et de Math\'ematiques Sociales\\
  CNRS, UMR 8557\\
  54 Bd Raspail, 75006 Paris\\
  France} \email{pom@ehess.fr}
\begin{abstract}
Classes of graphs with {\em bounded expansion} are a
generalization of both proper minor closed classes and degree bounded
classes. Such classes are based on a new invariant, the
{\em greatest reduced average density (grad) of $G$ with rank $r$},
$\rdens{r}(G)$. These classes are also characterized by the existence of several
partition results such as the existence of low tree-width and low
tree-depth colorings \cite{ICGT05}\cite{POMNI}. 
These results lead to several new linear time
algorithms, such as an algorithm for counting all the isomorphs of a
fixed graph in an 
input graph or an algorithm for checking whether there exists a subset
of vertices of 
{\em a priori} bounded size such that the subgraph induced by this
subset satisfies some arbirtrary but fixed first order sentence. We
also show that for fixed $p$, computing the distances between two
vertices up to distance $p$ may be performed in constant time per
query after a linear time preprocessing. We also show, extending
several earlier results, that a class of graphs has
sublinear separators if it has sub-exponential expansion. This result
result is best possible in general.
\end{abstract}

\maketitle
\section{Introduction}
The concept of tree-width \cite{Halin},\cite{GM1},\cite{Wagner} is
central to the analysis of graphs with forbidden minors done by
Robertson and Seymour and gained much algorithmic attention
thanks to the general complexity result of Courcelle about
monadic second-order logic graph properties decidability for graphs
with bounded tree-width \cite{Courcelle1},\cite{Courcelle2}. It
appeared that many NP-complete problems may be solved in polynomial
time when restricted to a class with bounded tree-width.
This restriction of tree-width is quite a strong one, as it 
does not include the class of planar graphs, for instance. 

Another way is to consider partitions of graphs into parts such that
any $p$ of them induce a graph with low tree-width. DeVos et
al. \cite{2tw} proved that for any proper minor closed class of graphs
$\mathcal C$ ---
that is: any minor closed class of graphs excluding at least one minor
--- and any integer $p$, there exists a constant $N(\mathcal C,p)$
so that any graph $G\in\mathcal C$ has a vertex-partition into at most
$N(\mathcal
C,p)$ parts such that
any $i\leq p$ parts induce a graph of tree-width at most $(i-1)$.

It is then natural to ask whether the parts could be choosen 
even ``smaller'' or ``simple''.
This issue has been studied in \cite{Taxi_tdepth} where the authors
introduce the {\em tree-depth} $\depth(G)$ of a graph $G$ as the
minimum height of a 
rooted forest including the graph in its closure. This minor monotone
invariant is related to tree-width by
$\tw(G)+1\leq\depth(G)\leq\tw(G)\log n$, where $n$ is the order of
$G$. The class of graphs with bounded tree-depth appears to be
particularly small, as it includes only a bounded number of rigid
graphs (that is: graphs having no non-trivial automorphisms) and as it
excludes long paths (to compare with classes with bounded tree-width
which exclude big grids). The main result of \cite{Taxi_tdepth} is
that for any proper minor closed class of graphs $\mathcal C$ and any
integer $p$, there exists an integer $N'(\mathcal C,p)$ such that any
graph $G\in\mathcal C$ has a vertex-partition into at most
$N'(\mathcal C,p)$ parts such that any $i\leq p$ parts induce a graph
of tree-depth at most $i$. It is also proved in \cite{Taxi_tdepth} 
that the tree-depth is
the greatest graph invariant for which such a statement holds.

Our first proof \cite{Taxi_tdepth} of this decomposition result relied
in the result 
of DeVos et al. and thus indirectly to the Structural Theorem of
Robertson and Seymour \cite{GM16}.
However since then, we generalized these results \cite{ICGT05}\cite{POMNI}
 to classes
with bounded expansion (which may be seen as a generalization of both
proper minor closed classes and degree bounded classes). 
Our prrof is both more general and conceptually easier.
Even better: it leads to a linear time algorithm that we shall
describe here. Our main goal will be then to show that this algorithm
has a wide range of algorithmic applications.

Before we shall consider algorithmic consequences, 
we shall introduce bounded
expansion and related concepts in Section \ref{sec:grad}.

In Section \ref{sec:aug} we describe the augmentation process
which is the basis of the partition theorem and propose
a linear time algorithm for it.

\section{The grad of a graph and classes with bounded expansion}
\label{sec:grad}
 The {\em
  distance} $d(x,y)$ between two vertices $x$ and $y$ of a graph is the minimum
length of a path linking $x$ and $y$, or $\infty$ if $x$ and $y$ do not
belong to the same connected component.
The {\em radius} $\rho(G)$ of a connected graph $G$ is:
$\rho(G)=\min_{r\in V(G)}\max_{x\in V(G)} {\rm d}(r,x)$

\begin{definition}
Let $G$ be a graph. A {\em ball} of $G$ is a subset of vertices
inducing a connected subgraph. 
The set of all the families of balls of $G$ is noted
$\mathfrak{B}(G)$.
The set of all the families of balls of $G$ including no two
intersecting balls is noted $\mathfrak{B}_1(G)$.

Let $\mathcal P=\{V_1,\dotsc,V_p\}$ be a family of balls of $G$.
\begin{itemize}
\item
 The {\em radius} $\rho(\mathcal P)$ of $\mathcal P$ is
 $\rho(\mathcal P)=\max_{X\in \mathcal P}\rho(G[X])$ 
\item
The {\em quotient} $G/\mathcal P$ of $G$ by 
  $\mathcal P$ is a graph with vertex set $\{1,\dotsc,p\}$ and edge
  set $E(G/\mathcal P)=\{\{i,j\}: (V_i\times V_j)\cap
  E(G)\neq\emptyset\text{ or }V_i\cap V_j\neq\emptyset\}$.
\end{itemize}
\end{definition}

\begin{definition} 
The {\em greatest reduced average density (grad) of $G$
      with rank $r$} is
$$\rdens{r}(G)=\max_{\substack{\mathcal P\in\mathfrak{B}_1(G)\\
\rho(\mathcal P)\leq r}}\frac{|E(G/\mathcal P)|}{|\mathcal P|}$$
\end{definition}

The first grad, $\rdens{0}$ , is closely related to
degeneracy ($G$ is $k$-degenerated iff $k\geq
\lfloor 2\rdens{0}(G)\rfloor$).  The grads of a graph form an
non decreasing sequence which becomes constant starting from some index (smaller than
the order of the graph).

\begin{definition}
A class of graphs $\mathcal C$ has {\em bounded expansion} if there
exists a function $f:\bbbn\rightarrow\bbbr$ such that for every graph $G \in \mathcal C$ and every $r$ holds
\begin{equation}
 \rdens{r}(G)\leq f(r)
\end{equation}
\end{definition}

Here are some examples of class with bounded expansion:
\begin{example}
Any proper minor closed class of graphs has expansion bounded by a
constant function. Conversely, any class of graphs with expansion
bounded by a constant is included in some proper minor closed class of
graphs.
\end{example}
\begin{proof}
If $\mathcal C$ is a proper minor closed class of graph, the graphs in
$\mathcal C$ are $k$-degenerated for some integer $k$ hence
$\rdens{r}(G)\leq k+1$ for any $G\in\mathcal C$.

Conversely, assume $\mathcal C$ is a class of graph with expansion
bounded by a constant $C$. Let $\mathcal C'$ be the class defined by
$\mathcal C'=\{G: \forall r\geq 0, \rdens{r}(G)\leq C\}$. This class
obviously includes $\mathcal C$. Let $G\in\mathcal C'$ and let $H$ be
a minor of $G$. Then for any $r\geq 0$,
$\rdens{r}(H)\leq\rdens{\card{V(G)}}(G)\leq C$ thus $H\in\mathcal C'$. Hence
$\mathcal C'$ is a proper minor closed class as it does not include $K_{2C+2}$ (as
$\rdens{0}(K_{2C+2})=C+1$).
\end{proof}

\begin{example}
Let $\Delta$ be an integer. Then the class of graphs with maximum
degree at most $\Delta$ has expansion bounded by the exponential function
$f(r)=\Delta^{r+1}$.
\end{example}

\begin{example}
In \cite{miller98} is introduced a class of graphs which occurs
naturally in finite-element and finite-difference problems. These graphs
correspond to graphs embedded in $d$-dimensional space in a certain
manner. It is proved in \cite{teng98} that these graphs excludes $K_h$
as a depth $L$ minor if $h=\Omega(L^d)$. Hence they form (for each
$d$) a class with polynomialy bounded expansion.
\end{example}

The next example show that the bounded function can be any arbitrary
increasing function:

\begin{example}
Let $f$ be any increasing function from $\bbbn$ to $\bbbn\setminus\{0,1,2\}$.
Then there exists a class $\mathcal C$ such that $\mathcal C$ has
expansion bounded by $f$ but by no smaller integral function.
\end{example}
\begin{proof}
Consider the class $\mathcal C$ whose elements are $K_4$ and the graphs $G_n$ obtained
by subdividing $2n$ times the complete graph $K_{2f(n)+1}$ (for $n\geq 1$).
As $2\leq \rdens{r}(G_n)<3$ for $r<n$ and as $\rdens{r}(G_n)=f(n)$ for
$r\geq n$, we conclude.
\end{proof}

\begin{example}
If $\mathcal C$ is a class with bounded expansion and
if $c$ is any fixed integer then the class $\mathcal C'\bullet K_c$
whose elements are the lexicographic products $G\bullet K_c,
G\in\mathcal C$ still has bounded expansion  \cite{POMNI}.

 It should be noted that
such a statement is false for proper minor closed classes in a strong
sense: for any $n\in\bbbn$, $K_n$ is a minor of ${\rm
Grid}(2n,2n)\bullet K_2$.
\end{example}

\subsection{Few properties of tree-depth}

A {\em rooted forest} is a disjoint union of rooted trees.
The {\em height} of a vertex $x$ in a rooted
forest $F$ is the number of vertices of a
path from the root (of the tree to which $x$ belongs to) to $x$ and is noted $\height(x,F)$.
The {\em height} of $F$ is the maximum height of the vertices of $F$.
Let $x,y$ be vertices of $F$. The vertex $x$ is an {\em ancestor} of $y$ in $F$ if $x$ belongs to
the path linking $y$ and the root of the tree of $F$ to which $y$ belongs to.
The {\em closure} $\clos(F)$ of a rooted forest $F$ is the graph with
vertex set $V(F)$ and edge set $\{\{x,y\}: x\text{ is an ancestor of
  }y\text{ in }F, x \neq y\}$. A rooted forest $F$ defines a partial order on its set of vertices:
$x\leq_F y$ if $x$ is an ancestor of $y$ in $F$.
The comparability graph of this partial order is obviously $\clos(F)$.

\begin{definition}
The {\em tree-depth} $\depth(G)$ of a graph $G$ is the minimum height
of a rooted forest $F$ such that $G\subseteq\clos(F)$.
\end{definition}

\begin{lemma}
\label{lem:fin}
Let $G$ be a connected graph with maximum degree $\Delta$ and
tree-depth $t\geq 1$.
Then $G$ has order $n\leq 1+\Delta+\dots+\Delta^{t-1}$.
\end{lemma}
\begin{proof}
We proceed by induction over $t$. If $t=1$, $G=K_1$ and
$\depth(K_1)=1$. Assume the inequality has been proved for graphs with
tree-depth at most $(t-1)$ with $t\geq 2$ and let $G$ be a connected
graph with tree-depth $t$.
As $G$ is connected it contains a vertex $r$ such that
$\depth(G-r)=\depth(G)-1=t-1$. Let $G_1,\dots,G_k$ be the connected
components of $G-r$. All of these have tree-depth at most $t-1$. By
induction, they have order at most $1+\Delta+\dots+\Delta^{t-2}$.
As $k\leq\Delta$, we conclude.
\end{proof}

\begin{lemma}
\label{lem:td_path}
For $k\geq 1, \depth(P_k)=\lceil\log_2 (k+1)\rceil$
\end{lemma}
\begin{proof}
According to 
lemma~\ref{lem:fin} a path of
tree-depth $t$ has order at most $1+2+\dots+2^{t-1}=2^t-1$. It follows that
the tree-depth of $P_k$ is at least $\log_2(k+1)$.

Moreover let $x_1,\dots,x_{2^t-1}$ be the vertices of a path of order
$2^t-1$ in the order in which they appear on the path. Let $w(i)$
be the base $2$ word of length $t$ corresponding to the number $i$ 
(for instance, if $t=3$, $w(1)=001, w(2)=010, \dots, w(7)=111$).
Let $c(i)$ be the rank of the rightmost $1$ of $w(i)$ (that is, for $t=3$,
$c(1)=c(3)=c(5)=c(7)=3$, $c(2)=c(6)=2$ and $c(4)=1$). Then $c$ is a
centered coloring of $P_{2^t-1}$ with $t$ colors thus
$\depth(P_{2^t-1})\leq t$. Finally we note that
 $\depth(P_n)$ increases with $n$.
\end{proof}

\begin{lemma}
\label{lem:path}
Let $G$ be a graph and let $P_k$ be the longest path in
$G$.

Then $\lceil\log_2 (k+1)\rceil\leq \depth(G)\leq \binom{k+2}{2}-1$.
\end{lemma}
\begin{proof}
 As the tree-depth is minor monotone, any graph including a path $P_k$
as a subgraph as 
tree-depth at least $\depth(P_k)=\lceil\log_2 (k+1)\rceil$ (according
to Lemma \ref{lem:td_path}).

Conversely, let us prove by induction over $k\geq 1$ that a graph which
includes no path $P_k$ has tree depth at most
$\binom{k+1}{2}$. Obviously the statement holds for $k=1$
(graphs without edges has tree-depth $1$). Assume the statement has
been proved up to $(k-1)$ for some $k\geq 2$. Let $G$ be a graph with
no path $P_k$. Without loss of generality we may assume
that $G$ includes a $P_{k-1}$ and that $G$ is connected
(as the tree-depth of a non-connected graph is the maximum of the
tree-depths of its connected components). Let $P$ be such a path of
$G$. Assume $G-V(P)$ includes some path $P'$ isomorphic to $P_{k-1}$. According
to the connectivity of $G$, there exists some minimum length path
$P''$ linking a vertex of $P$ to a vertex of $P'$ (and this path has
length at least $1$). Then $P\cup P'\cup P''$ includes a $P_k$, a
contradiction. 
Thus $G-V(P)$ includes no
$P_{k-1}$. By induction, $\depth(G-V(P))\leq\binom{k}{2}-2$
hence $\depth(G)\leq \depth(G-V(P))+\card{V(P)}\leq
\binom{k}{2}+k=\binom{k+1}{2}$. If follows that if $P_k$ is the longest
path in $G$, $\depth(G)$ is at most $\binom{k+2}{2}-1$.
\end{proof}

\section{Basics}
We shall first mention some basic linear time
algorithms, as well as the basic data structures used for input and
output of our algorithms. Concerning the data structure used for the
computations, any standard one will do, but we will have in
mind the simple data structure of {\tt PIGALE} library
\cite{Taxi_pigale}. 

\subsection{Digraph representation}

A computed directed graph $\vG$ will be represented as an array $D$ of lists
indexed by integers $1,\dotsc,n$. In the list $D[i]$ will be gathered
all the couples $(j,e)$  such that $(j,i)$ is an arc of $\vG$ 
with index $e\in\{1,\dotsc,m\}$ (where $m$ is the size of $\vG$). 
This representation can be easily constructed from any standard one
in linear time. Moreover, it is possible to filter out parallel edges in linear
time using bucket-sort, and to transform into any standard representation in
linear time. The main interest in numbering the vertices and edges stands in the 
possibility to use ``raw'' integer arrays to store any needed information
and to ease bucket-sorting (this simple fact is central to the efficiency of
Pigale's data structure \cite{Taxi_pigale}).

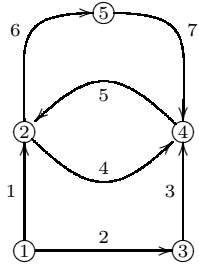
$$
\begin{xy}<15pt,0cm>:
(0,0)*+{\scriptstyle 1}*[o]{\cir{}}="va",
(0,3)*+{\scriptstyle 2}*[o]{\cir{}}="vb",
(4,0)*+{\scriptstyle 3}*[o]{\cir{}}="vc",
(4,3)*+{\scriptstyle 4}*[o]{\cir{}}="vd",
(2,6)*+{\scriptstyle 5}*[o]{\cir{}}="ve"
\ar@{->}^{1} "va";"vb"
\ar@{->}^{2} "va";"vc"
\ar@{->}^{3} "vc";"vd"
\ar@{->}^{4}@(dr,dl) "vb";"vd"
\ar@{->}^{5}@(ul,ur) "vd";"vb"
\ar@{->}^{6}@(u,l) "vb";"ve"
\ar@{->}^{7}@(r,u) "ve";"vd"
\end{xy}
$$
$$
\begin{aligned}
D[1]&=()\\
D[2]&=((1,1),(4,5))\\
D[3]&=((1,2))\\
D[4]&=((2,4),(3,3),(5,7))\\
D[5]&=((2,6))
\end{aligned}
$$

Notice that the used representation of a directed graph $\vG$ allows to
answer the question ``is there an arc from vertex $i$ to vertex $j$''
in time $O(\md(\vG))$. Notice that this simple observation has by
itself many algorithmic consequences \cite{chrobak}.

\subsection{Low indegree orientation}

The aim of the following algorithm is to compute a low-indegree
orientation of the graph with vertex set $\{1,\dots,n\}$ and list of
edges $L$.

\begin{lemma}
Let $G$ be a graph of order $n$ and size $m$. There is an
$O(n+m)$-time algorithm which computes an acyclic orientation of $G$
with maximum indegree $\lfloor 2\rdens{0}(G)\rfloor$.
\end{lemma}
\begin{proof}
First we 
compute a representation of the graph in any suitable data structure
like {\tt PIGALE}'s data structure \cite{Taxi_pigale}. All of this may
be easily done in time $O(m)$. Then we do the following:

\vspace{5mm}
\begin{algorithmic}
\ENSURE $D$ represents an orientation $\vG$ of $G$ such that
$\md(\vG)\leq\lfloor 2\rdens{0}(G)\rfloor$. 
\STATE Let $D[1\dots n]\leftarrow ()$.
\STATE $\forall v: d[v]\leftarrow$ degree of $v$ in $G$
\STATE $\forall i: T[i]\leftarrow$ list of vertices of $G$ with degree $i$
\STATE $\delta\leftarrow 0$
\STATE $m\leftarrow 0$
\WHILE {$\delta<n$}
  \IF {$T[\delta]\ne ()$}
    \STATE pop $v$ out of $T[\delta]$.
    \STATE let $d[v]\leftarrow 0$
      \FORALL {$w$ neighbour of $v$}
         \IF{ $d[w]>0$}
            \IF {$d[w]>\delta$}
               \STATE extract $w$ from  $T[d[w]]$
               \STATE insert $w$ in $T[d[w]-1]$
            \ENDIF
            \STATE let $d[w]\leftarrow d[w]-1$
            \STATE $m\leftarrow m+1$
            \STATE append $(w,m)$ to $D[v]$
         \ENDIF
      \ENDFOR
  \ELSE
    \STATE $\delta\leftarrow\delta+1$
  \ENDIF
\ENDWHILE
\end{algorithmic}

In this algorithm, if $\delta$ is increased the the 
subgraph of $G$ induced by the remaining vertices
has minimum degree greater than $\delta$. It follows
that the maximum value of $\delta$ reached by the algorithm is less or
equal to the maximum average degree of $G$, that is: $\delta\leq
2\rdens{0}(G)$.
It follows that this algorithm computes an acyclic orientation of $G$
with maximum indegree $\lfloor 2\rdens{0}(G)\rfloor$ in time $O(m)$.
\end{proof}

\section{Transitive fraternal augmentations of graphs in linear time}
\label{sec:aug}
\subsection{Theory}
In the following, a directed graph $\vG$ may not have a loop and for
any two of its vertices $x$ and $y$, $\vG$ includes at most one arc
from $x$ to $y$ and at most one arc from $y$ to $x$ 
(thus at most two arcs may connect $x$ and $y$, one in each direction).

\begin{definition}
Let $\vG$ be a directed graph. A {\em $1$-transitive fraternal augmentation} of $\vG$
is a directed graph $\vH$ with the same vertex set, including all the
arcs of $\vG$ and such that, for any distinct vertices $x,y,z$,
\begin{itemize}
\item if $(x,z)$ and $(z,y)$ are arcs of $\vG$ then $(x,y)$ is an arc
  of $\vH$ ({\em transitivity}),
\item  if $(x,z)$ and $(y,z)$ are arcs of $\vG$ then $(x,y)$ or $(y,x)$ is an arc
  of $\vH$ ({\em fraternity}).
\end{itemize}

A {\em transitive fraternal augmentation} of a directed graph $\vG$ is a sequence
$\vG=\vG_1\subseteq\vG_2\subseteq\dotsb\subseteq \vG_i\subseteq\vG_{i+1}\subseteq\dotsb$,
such that $\vG_{i+1}$ is a $1$-transitive fraternal augmentation of $\vG_i$ for any
$i\geq 1$.
\end{definition}

The key result of \cite{POMNI} claims the existence of density bounded
transitive fraternal augmentations:
\begin{lemma}[Special case of Lemma 6.1 of \cite{POMNI}]
\label{lem:faug}
There exists polynomials $P_i\ (i\geq 0)$ such that for any directed graph
 $\vG$ and any $1$-transitive fraternal augmentation $\vH$ 
of $\vG$ we have
\begin{equation}
\rdens{r}(H) \leq P_{2r+1}(\md(\vG)+1,\rdens{2r+1}(G)),
\end{equation}
where $G$ and $H$ stand for the simple undirected graphs underlying $\vG$
and $\vH$.
\end{lemma}

Although quite technical, the next result is a simple direct
consequence of Lemma~\ref{lem:faug}:

\begin{corollary}
Let $\mathcal C$ be a class with expansion bounded by a function $f$
and let $F:\bbbn^2\rightarrow\bbbn$. 

Define $A(r,i)$ and $B(i)$ recursively as follows (for $i\geq 1$
and $r\geq 0$):
\begin{align*}
A(r,1)&=f(r)\\
B(1)&=2f(0)\\
A(r,i+1)&=P_{2r+1}(B(i)+1,A(2r+1,i))\\
B(i+1)&=F(B(i),A(0,i+1))
\end{align*}

Assume $G\in\mathcal C$ and $\vG=\vG_1\subseteq\vG_2\subseteq\dotsb\subseteq
\vG_i\subseteq\vG_{i+1}\subseteq\dotsb$ is a transitive fraternal
augmentation of $G$ such that $\md(\vG_{i+1})\leq F(\md(\vG_i),\rdens{0}(G_{i+1}))$
(for $i\geq 1$) and such that $\md(\vG_1)\leq 2f(0)$.
Then:
\begin{align*}
\rdens{r}(G_i)&\leq A(r,i)\\
\md(\vG_i)&\leq B(i)
\end{align*}
\end{corollary}

We now present a linear time implementation of this procedure, where
 it will be checked that $\md(\vG_{i+1})\leq
 \md(\vG_i)^2+2\rdens{0}(G_i)$, that is: $F(x,y)=x^2+2y$.

\subsection{The algorithm for one step augmentation}
In the augmentation process, we add two kind of arcs:
transitivity arcs and fraternity arcs. Let us start with transitivity
ones:

$$
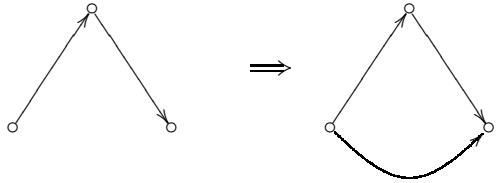
\begin{xy}<15pt,0cm>:
(0,0)*[o]{\circ{}}="va",
(4,0)*[o]{\circ{}}="vb",
(2,3)*[o]{\circ{}}="vc",
(6,1.5)="a",
(7,1.5)="b",
(8,0)*[o]{\circ{}}="wa",
(12,0)*[o]{\circ{}}="wb",
(10,3)*[o]{\circ{}}="wc",
\ar@{=>} "a";"b"
\ar@{->} "va";"vc"
\ar@{->} "vc";"vb"
\ar@{->} "wa";"wc"
\ar@{->} "wc";"wb"
\ar@{->}@(dr,dl) "wa";"wb" 
\end{xy}
$$

\vspace{5mm}
\begin{algorithmic}
\REQUIRE $D$ represents the directed graph to be augmented.
\ENSURE $D'$ represents the array of the added arcs.
\STATE Initialize $D'$.
\FORALL{$v\in\{1,\dots,n\}$}
\FORALL{$(u,e)\in D[v]$}
\FORALL{$(x,f)\in D[u]$}
\STATE $m\leftarrow m+1$; append $(x,m)$ to $D'[v]$.
\ENDFOR
\ENDFOR
\ENDFOR
\end{algorithmic}
\vspace{5mm}

This algorithm runs in $O(\md(\vG)^2n)$ time, where $\md(\vG)$ is the
maximum indegree of the graph to be augmented. It computes the list
array $D'$ of
the transitivity arcs which are missing in $\vG$, missing arcs
may appear more than once in the list, but the number of added edges
cannot exceed $\md(\vG)^2n$.

Now, we shall consider the fraternity edges.

$$
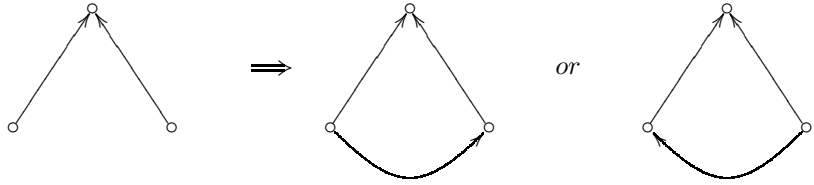
\begin{xy}<15pt,0cm>:
(0,0)*[o]{\circ{}}="va",
(4,0)*[o]{\circ{}}="vb",
(2,3)*[o]{\circ{}}="vc",
(6,1.5)="a",
(7,1.5)="b",
(8,0)*[o]{\circ{}}="wa",
(12,0)*[o]{\circ{}}="wb",
(10,3)*[o]{\circ{}}="wc",
(14,1.5)*+{or},
(16,0)*[o]{\circ{}}="za",
(20,0)*[o]{\circ{}}="zb",
(18,3)*[o]{\circ{}}="zc"
\ar@{=>} "a";"b"
\ar@{->} "va";"vc"
\ar@{->} "vb";"vc"
\ar@{->} "wa";"wc"
\ar@{->} "wb";"wc"
\ar@{->}@(dr,dl) "wa";"wb" 
\ar@{->} "za";"zc"
\ar@{->} "zb";"zc"
\ar@{->}@(dl,dr) "zb";"za" 
\end{xy}
$$

\vspace{5mm}
\begin{algorithmic}
\REQUIRE $D$ represents the directed graph to be augmented.
\ENSURE $L$ represents the list of edges to be added.
\STATE $L=()$.
\FORALL{$v\in\{1,\dots,n\}$}
\FORALL{$(x,e)\in D[v]$}
\FORALL{$(y,f)\in D[v]$}
\IF{$x<y$}
\STATE append $(x,y)$ to $L$.
\ENDIF
\ENDFOR
\ENDFOR
\ENDFOR
\end{algorithmic}
\vspace{5mm}

This algorithm runs in $O(\md(\vG)^2n)$-time and computes the list of
the fraternity edges, edges may appear more than once but the length of the
list $L$ cannot exceed $\md(\vG)^2n/2$.

The simplification of $L$, the computation of a low indegree orientation
of the edges in $L$ and the merge/simplification with the arcs in $D$ and $D'$
may be achieved in linear time (precisely: in $O(\md(\vG)^2 n)$-time).
\begin{theorem}
\label{th:augmlin}
For any class $\mathcal C$ with bounded expansion and any fixed 
integer $c$,
there exists an algorithm which computes, given an input graph
$G\in\mathcal C$, a
transitive fraternal augmentation
$\vG=\vG_1\subseteq\vG_2\subseteq\dotsb\subseteq \vG_c$ of $G$ in time
$O(n)$.
\end{theorem}

\section{Distances}
The following result is a weighted extention of
the basic observation that bounded orientations allows $O(1)$-time
checking of adjacency \cite{chrobak}.

\begin{theorem}
For any class $\mathcal C$ with bounded expansion and for any integer
$k$, there exists a linear time preprocessing algorithm so that for
any preprocessed
$G\in\mathcal C$ and any pair $\{x,y\}$ of vertices of $G$ the value
$\min(k,{\rm dist}(x,y))$ may be computed in $O(1)$-time.
\end{theorem}
\begin{proof}
The proof goes by a variation of our augmentation algorithm so that
each arc $e$ gets a weight $w(e)$ and each added arc gets weight
$\min(w(e_1)+w(e_2))$ over all the pairs $(e_1,e_2)$ of arcs which may
imply the addition of $e$ and 
simplification should keep the minimum weighted arc.

Then, after $k$ augmentation steps, two vertices at distance at most
$k$ have distance at most $2$ in the augmented graph. The value
$\min(k,{\rm dist}(x,y))$ then equals
$\min(k,w((x,y)),w((y,x)),\min_{(z,x),(z,y)\in\vG}(w(z,x)+w(z,y)))$.
\end{proof}

\section{$p$-centered colorings and tree-decomposition}
\subsection{Theory}

\begin{definition}
A {\em tree-decomposition} of a graph $G$ consists in a pair $(T,\lambda)$
formed by a tree $T$ and a function $\lambda$ mapping vertices of $T$
to subsets of $V(G)$ so that for all $v\in V(G), \{x\in V(T):
v\in\lambda(x)\}$ induces a subtree of $T$, and such that
for any edge $\{v,w\}$ of $G$ there exists $x\in V(T)$ such that
$\{v,w\}\subseteq\lambda(x)$. 

The {\em width} of a tree decomposition  $(T,\lambda)$ is $\max_{v\in
V(G)}\card{\lambda(v)}-1$. The {\em tree-width} of $G$ is the minimum
width of any tree-decomposition of $G$.
\end{definition}

From a rooted tree $Y$ of height at most $p$ such that
$G\subseteq\clos(Y)$ it is straightforward to construct a
tree-decomposition $(T,\lambda)$ of $G$ having width at most $(p-1)$:
Set $T=Y$ and define $\lambda(x)=\{v\leq _Y x\}$. Then for any $v$,
$\{x\in V(T): v\in\lambda(x)\}=\{x\geq_Y v\}$ induces the subtree of
$Y$ rooted at $v$ (hence a subtree of $T$). Moreover, as $G\subseteq\clos(Y)$, any edge
$\{x,y\}$ with $x<_Y y$ is a subset of $\lambda(y)$. Hence
$(T,\lambda)$ is a tree-decomposition of $G$. As $\max_{v\in
V(G)}\card{\lambda(v)}=\height(Y)\leq p$, this tree-decomposition has
width at most $(p-1)$. Last, this tree-decomposition may be obviously
constructed in linear time.

\begin{definition}
A {\em centered coloring} of a graph $G$ is a coloring of the vertices
such that in any connected subgraph some color appears exactly once.

For an integer $p$, a {\em $p$-centered coloring} of $G$ is a coloring
of the vertices such that in any connected subgraph either some color
appears exactly once, or at least $p$ different colors appear.
\end{definition}

\subsection{The algorithm}

\vspace{5mm}
\begin{algorithmic}
\REQUIRE $c$ is a centered-coloring of the graph $G$ using colors
$1,\dotsc,p$.
\ENSURE $\mathcal F$ is a rooted forest such that
$G\subset\clos(\mathcal F)$.
\STATE Set $\mathcal F=\emptyset$.
\STATE Let ${\rm Big}[\ ]$ be an array of size $p$.
\FORALL{Connected component $G_i$ of $G$}
\STATE Initialize  ${\rm Big}[\ ]$ to {\rm\bf false}.
\STATE Set ${\rm root\_color}\leftarrow 0$.
\FORALL{$v\in V(G_i)$}
\IF{${\rm Big}[c[v]]={\rm\bf false}$}
\IF{$c[v]={\rm root\_color}$}
\STATE ${\rm  root\_color}\leftarrow 0, {\rm
Big}[c[v]]\leftarrow {\rm\bf true}$.
\ELSE 
\STATE ${\rm root}\leftarrow v; {\rm  root\_color}\leftarrow c[v]$.
\ENDIF
\ENDIF
\ENDFOR
\STATE Recurse on $G-{\rm root}$ thus getting some rooted forest $\mathcal
F'=\{Y_1',\dots,Y_j'\}$.
\STATE Add to $\mathcal F$ the tree with root ${\rm root}$ and subtrees
$Y_1,\dots,Y_j$, where the sons of ${\rm root}$ are the roots of
$Y_1,\dots,Y_j$.
\ENDFOR
\end{algorithmic}
\vspace{5mm}

This algorithms clearly runs in $O(pm)$ time. If $G$ is connected, it
returns a rooted tree $Y$ of height at most $p$ such that $G\subseteq\clos(Y)$.

\section{Application to subgraph isomorphism problem}
For general subgraph isomorphism problem of deciding wether a graph
$G$ contains a subgraph isomorphic to a graph $H$ of order $l$, 
the better known general bound is $O(n^{\alpha l/3})$ where $\alpha$ is
the exponent of square matrix fast multiplication algorithm
\cite{NP85} (hence $O(n^{0.792\ l})$ using the fast matrix algorithm of \cite{coppersmith90}).
 The particular case of
subgraph isomorphism in planar graphs have been studied by Plehn and
Voigt \cite{plehn91}, Alon \cite{alon95} with superlinear bounds and
then by Eppstein \cite{Epp-SODA-95}\cite{Epp-JGAA-99} who gave the
first linear time algorithm for fixed pattern $H$ and $G$ planar and
then extended his result to graphs with bounded genus \cite{Epp-Algo-00}.
We generalize this to classes with bounded expansion.

We shall now make use of the following result for graphs with bounded
tree-width: 

\begin{lemma}[Eppstein, Lemma 2 of \cite{Epp-JGAA-99}]
\label{lem:Epp2}
Assume we are given graph $G$ with $n$ vertices along with a
tree-decomposition $T$ of $G$ with width $w$. Let $S$ be a subset of
vertices of $G$, and let $H$ be a fixed graph with at most $w$
vertices. Then in time $2^{O(w\log w)}n$ we can count all isomorphs of
$H$ in $G$ that include some vertex in $S$. We can list all such
isomorphs in time  $2^{O(w\log w)}n+O(kw)$, where $k$ denotes the
number of isomorphs and the term $kw$ represents the total output
size.
\end{lemma}

We shall prove here the following extention of the results of
\cite{Epp-JGAA-99}\cite{Epp-Algo-00}: 

\begin{theorem}
Let $\mathcal C$ be a class with bounded expansion and let $H$ be a
fixed graph. Then there exists a linear time algorithm
which computes, from a pair $(G,S)$ formed by a graph $G\in\mathcal C$
and a subset $S$ of vertices of $G$, the number of isomorphs of $H$ in
$G$ that include some vertex in $S$. There also exists an algorithm
running in time $O(n)+O(k)$ listing all such isomorphism where  $k$ denotes the
number of isomorphs (thus represents the output size).
\end{theorem}
\begin{proof}
This is a direct consequence of Theorem~\ref{th:augmlin} and
Lemma~\ref{lem:Epp2}.
\end{proof}
\section{Local decidability problems}

Monadic second-order logic (MSOL) is an extention of first-order logic (FOL)
that includes vertex and edge sets and belonging to these
sets. The following theorem of Courcelle has been applied to solve
many optimization problems. 

\begin{theorem}[Courcelle \cite{Courcelle1}\cite{Courcelle2}]
\label{th:msol}
Let $\mathcal K$ be class of finite graphs $G=\langle V,E,R\rangle$
represented as $\tau_2$-structures, that is: by two sorts of
elements (vertices $V$ and edges $E$) and an incidence relation $R$,
and $\phi$ be a MSOL($\tau_2$) sentence. If $\mathcal K$ has bounded
tree width and $G\in\mathcal K$, then checking wether $G\vDash\phi$
can be done in linear time.
\end{theorem}

Combining Theorem \ref{th:msol} with Theorem \ref{th:augmlin}, we get:

\begin{theorem}
Let $\mathcal C$ be a class with bounded expansion and let $p$ be a
fixed integer. Let $\phi$ be a FOL($\tau_2$) sentence. Then there exists a
linear time algorithms to check  $\exists X: (\card{X}\leq p)\wedge(G[X]\vDash\phi)$.
\end{theorem}

Thus for instance:

\begin{theorem}
Let $\mathcal K$ be a class with bounded expansion and let $H$ be a
fixed graph. Then, for each of the next properties there exists a linear
time algorithm to decide whether a graph $G\in\mathcal K$ satisfies them:
\begin{itemize}
\item $H$ has a homomorphism to $G$,
\item $H$ is a subgraph of $G$,
\item $H$ is an induced subgraph of $G$.
\end{itemize}
\end{theorem}

Although there is an (easy) polynomial algorithm to decide whether
$td(G) \leq k$ for any fixed $k$, if P$\neq$NP then no
polynomial time approximation algorithm for the tree-depth can
guarantee an error bounded by $n^\epsilon$, where $\epsilon$ is a
constant with $0<\epsilon<1$ and $n$ is the order of the graph
\cite{treewidth_approx}.
We shall now prove that the decision problem $td(G) \leq k$ for any
fixed $k$ may actually be decided in linear time:

\begin{lemma}
\label{lem:DFS}
Any Depth-First Search (DFS) tree $Y$ of connected graph $G$
satisfies:
\begin{itemize}
\item $G\subseteq\clos(Y)$,
\item $\depth(G)\leq\height(Y)\leq 2^{\depth(G)}-1$.
\end{itemize}
\end{lemma}
\begin{proof}
According to the basic properties of the DFS, a vertex $v$ of $G$ may not be
adjacent in $G$ to a vertex which is not comparable to $v$ with
respect to the tree order induced by the DFS tree $Y$ thus
$G\subseteq\clos(Y)$ and $\depth(G)\leq\height(Y)$.
 Moreover, $G$ includes $P_{\height(Y)}$ as a
subgraph (take any maximal tree chain) thus $\height(Y)\leq
2^{\depth(P_{\height(Y)})}-1$, according to Lemma \ref{lem:td_path}.
Hence
$\height(Y)\leq 2^{\depth(G)}-1$ as
$\depth(P_{\height(Y)})\leq\depth(G)$.
\end{proof}
\begin{theorem}
For any fixed $k$, there exists a linear time algorithm which decides
wether an input graph $G$ has tree-depth at most $k$ or not.
\end{theorem}
\begin{proof}
Without loss of generality we may assume $G$ is connected (for
otherwise we process all the connected components one by one).
Any DFS tree $Y$ of $G$ may be computed in $O(m)$ time, where $m$ is
the size of $G$. If $\height(Y)\geq 2^k$, the answer is ``No''
according to Lemma \ref{lem:DFS}.
Otherwise, consider the following
sentence $\Phi$:
$$
\begin{aligned}
\exists V_1\exists V_2\dots\exists V_k: &(\forall x\in V_1\ \forall y\in V_2, x\neq y)\wedge\dots\\
&\wedge (\forall x (\exists y\in V_1, x=y)\vee \dots )\\
&\begin{aligned}
\wedge (\forall A (\exists B &(\forall x\in A\ (x\in B))\wedge (\forall
x\in B\ 
\forall y\in A\ (y\in B)\vee\urcorner {\rm Adj}(x,y)))\\
&\vee (\exists x\in V_1\ (x\in A) \wedge (\forall y\in A\
(x=y)\vee\urcorner(y\in V_1)))\\
&\vee\dots\\
&\vee (\exists x\in V_k\ (x\in A) \wedge (\forall y\in A\
(x=y)\vee\urcorner(y\in V_1))))
\end{aligned}
\end{aligned}
$$

The first two lines express that $V_1,\dots,V_k$ shall be a partition
of the vertex set, and the next ones express that for any subset $A$
of vertices, either $G[A]$ is not connected of for some $i$ $A$
includes exactly one element of $V_i$, that is: $V_1,\dots,V_k$ is a
centered coloring of $G$. Such a centered coloring with $k$ colors
exists if and only if $G$ has tree-depth has most $k$
\cite{Taxi_tdepth}. It follows that $G\vDash\Phi$ if and only if
$\depth(G)\leq k$. As we only check $\Phi$ on graphs with tree depth
at most $2^k$ (given togather with a tree-decomposition easily deduced
from the DFS tree) and as $\Phi$ obviously belongs to $MSOL$, there exists,
according to Theorem \ref{th:msol}, a linear time algorithm to check
wether $G$ satisfies $\Phi$.
\end{proof}

\section{Vertex separators}
A celebrated theorem of Lipton and Tarjan \cite{lipton79} states that
any planar graph has a separator of size $O(\sqrt{n})$. Alon, Seymour
and Thomas \cite{AST90} showed that excluding $K_h$ as a minor ensures
the existence of a separator of size at most
$O(h^{3/2}\sqrt{n})$. Gilbert, Hutchinson, and Tarjan \cite{gilbert84} further proved
that graphs with genus $g$ have a separator of size $O(\sqrt{gn})$
(this result is optimal).
Plotkin et al. \cite{shallow} introduced the concept of {\em
limited-depth minor} exclusion and have shown that exclusion of small
limited-depth minors implies the existence of a small separator.
Precisely, they prove that any graph excluding $K_h$ as a depth $l$
minor has a separator of size $O(lh^2\log n+n/l)$ hence proving that
excluding a $K_h$ minor ensures the existence of a separator of size 
$O(h\sqrt{n}\log n)$. 

We use the following result to show that any class of graphs with 
sub-exponential
expansion has separators of sublinear size.

\begin{theorem}[Plotkin et al. \cite{shallow}]
\label{th:plot}
Given a graph with $m$ edges and $n$ nodes, and integers $l$ and $h$,
there is an $O(mn/l)$ time algorithm that will either produce a
$K_h$-minor of depth at most $l \log n$ or will find a separator of
size at most $O(n/l+4lh^2\log n)$.
\qed
\end{theorem}

\begin{lemma}
There exists a constant $C$ such that any 
graph $G$ has a separator of size at most $C\frac{n\log n}{z}$
whenever $z$ is an integer such that
\begin{equation}
2z(\rdens{z}(G)+2)\leq\sqrt{n\log n}.
\label{eq:z}
\end{equation}
\end{lemma}
\begin{proof}
Let $l=z/\log n$ and let $h=\lfloor \rdens{z}(G)+2\rfloor$. As
$\rdens{z}(G)\leq f(z)<h-1$, $G$ has no $K_h$ minor of depth at most
$l \log n$. According to Theorem~\ref{th:plot}, $G$ has a separator
of size at most $(C/2)(n/l+4lh^2\log n)$ for some fixed constant $C$, 
i.e. a separator of size at most
$(C/2)(\frac{n\log n}{z}+4z(\rdens{z}(G)+2)^2)\leq
C\frac{n\log n}{z}$.
\end{proof}
\begin{theorem}
\label{th:vs}
Let $\mathcal C$ be a class of graphs with expansion bounded
by a function $f$ such that $\log f(x)=o(x)$.

Then the graphs in $\mathcal C$ have separators of size $o(n)$.
\end{theorem}
\begin{proof}
Let $g(x)=\frac{\log f(x)}{x}$. By assumption,
$g(x)=o(1)$.
Define $\zeta(n)$ as the greatest integer such that
\begin{equation*}
\log f(\zeta(n))<\frac{\log n}{3}
\end{equation*}
Notice that $\zeta$ is increasing and $\lim_{n\rightarrow\infty}\zeta(n)=\infty$.
From the definition of $g(x)$, we deduce
 $\zeta(n)=\frac{\log f(\zeta(n))}{g(\zeta(n))}=\frac{\log
n}{3g(\zeta(n))}=o(\log n)$.
Thus $\log (2\zeta(n)(f(\zeta(n))+2))<\frac{\log
n}{3}(1+o(1))$. It follows that if $n$ is
sufficiently large (say $n>N$),  $\log (2\zeta(n)(f(\zeta(n))+2))<\frac{\log
n+\log\log n}{2}$, that is: $2\zeta(n)(f(\zeta(n))+2)<\sqrt{n\log
n}$. Thus if $n>N$, $G$ has a separator of size at most
$C\frac{n\log n}{\zeta(n)}=3g(\zeta(n))n=o(n)$.
\end{proof}

As random cubic graphs almost surely have bisection width at least
$0.101n$ (Kostochka and Melnikov, 1992), they have almost surely no
separator of size smaller than $n/20$ It follows that if $\log
f(x)=(\log 2)x$, the graphs have no sublinear separators any more.
This shows the optimality of Theorem \ref{th:vs}.
\begin{ack}
The authors would like to thank Martin Mares for his fruitful comments.
\end{ack}

\providecommand{\bysame}{\leavevmode\hbox to3em{\hrulefill}\thinspace}
\providecommand{\MR}{\relax\ifhmode\unskip\space\fi MR }
\providecommand{\MRhref}[2]{%
  \href{http://www.ams.org/mathscinet-getitem?mr=#1}{#2}
}
\providecommand{\href}[2]{#2}

\end{document}